\documentclass[10pt]{article}
\usepackage[utf8]{inputenc}
\usepackage[T1]{fontenc}
\usepackage{multirow}
\usepackage{amsmath}
\usepackage{amssymb}
\usepackage{soul}
\usepackage[mathscr]{eucal}
\usepackage{color}
\usepackage{authblk}

\newtheorem{Lemma}{Lemma}

\def\bx{{\bf x}}

\def\ve{{\varepsilon}}

\def\div{\mbox{{\rm div}}}

\def\va{\raise 2pt\hbox{,}}

\def\@xthm#1#2{\@beginassumption{#2}{\csname the#1\endcsname}{}\ignorespaces}
\def\@ythm#1#2[#3]{\@opargbeginassumption{#2}{\csname the#1\endcsname}{#3}\ignorespaces}%
\def\@beginassumption#1#2#3{\par\addvspace{8pt plus3pt minus2pt}%
              \noindent{\csname#1headfont\endcsname#1\ \ignorespaces#3 #2.}%
              \csname#1font\endcsname\hskip.5em\ignorespaces}
\def\@endassumption{\par\addvspace{8pt plus3pt minus2pt}\@endparenv}
\begin{document}


\title{Micro-Macro Derivation of Virus-Chemotaxis  Models}

\author[$\dagger$]{D. Burini}
\author[$\star$]{N. Chouhad}

\affil[$\dagger$]{University of Perugia, Perugia, Italy. \textit{dilettaburini@alice.it}}
\affil[$\star$]{Cadi Ayyad University, Ecole Nationale des Sciences Appliqu\'ees, Marrakech, Morocco. \textit{chouhadn@gmail.com}}

\maketitle

\begin{abstract}
This paper deals with the micro-macro derivation of  virus models  coupled with a reaction diffusion models that generates the dynamics in space  of the virus particles.   The first part  of the  presentation focuses, starting from~\cite{[BPTW19],[BT20]} on a survey and a critical analysis of some phenomenological models known in the literature. The second part shows how methods of the kinetic theory can be used to model the dynamics of the system  treated in our paper. The third part deals with the derivation of macroscopic models from the underlying description, delivered within a general framework of the kinetic theory.
\vskip.2cm
\noindent \textbf{keywords:} Kinetic theory, active particles, cross diffusion, multiscale methods.

\end{abstract}

\section{Aims and plan of the paper}\label{Sec:1}

The phenomenological derivation of models of biological tissues  can be obtained, at the macroscopic scale, by conservation, or equilibrium, equations that involve locally averaged quantities deemed to describe the state of the system.  These structures can be closed by means of heuristic models of the material behavior of the physical system under consideration. Different models correspond to each specific closure that are obtained by material  models generally valid only for physical conditions closed to equilibrium, while dynamical models are required to describe physical reality far from equilibrium.

Alternative methods have been  developed to tackle the aforementioned conceptual difficulty. Specifically, we refer to the method proposed in~\cite{[BC17],[BC19]} where the derivation of reaction-diffusion equations and of the celebrated Keller-Segel model~\cite{[KS70],[KS71]} was considered. This derivation  is somehow inspired to the  Hilbert's sixth problem, see\cite{[HILBERT]} which suggests the search of a unified approach to physical theories at all  representation scales, see also~\cite{[GK14]}.

In general,  the dynamics at the low scale is modeled by the kinetic theory of active particles~\cite{[BBDG21]}. This equation is expanded in terms of a small parameter corresponding to the mean distance between pair of particles, which is closed by neglecting the contribution of  higher order terms. The macro-scale model is obtained by taking low order terms of the expansion. This method is reviewed in the survey paper~\cite{[BBNS12]} mainly devoted to the micro-macro derivation Keller and Segel type models  (in short KS model).

In more detail, our paper is devoted to the micro-macro derivation of reaction-diffusion and cross-diffusion models in which a virus model is coupled with a reaction diffusion models that regulate the spatial dynamics of viral particles. This class of models can be defined \textit{exotic} by a term used to denote models in fields of behavioral sciences, for instance social and economical sciences, namely sciences where individual behaviors have an influence on the mechanical dynamics. Occasionally, the term \textit{model in complex environments} is used to denote the interaction between a \textit{first model}, essentially chemotaxis or cross diffusion, with a \textit{second (additional) model} which describes the external dynamics. The survey~\cite{[BOSTW21]} reports about a broad variety of this type of models.

In more details on the contents of our paper, Section 2 provides a description  and the phenomenological derivation of the aforementioned class of models. Section 3 presents the methodological approach used for the micro-macro derivation. Section 4 shows how the approach can be applied to the derivation of the macroscopic description of the models reported in Section 2. Finally, a critical analysis is presented in the last section looking ahead to research perspectives.

\section{Heuristic derivation of cross diffusion of virus models}\label{Sec.2}

The class of models presented in this section describes the dynamics of a virus model, where  space dynamics is induced by a transport mechanism which is modeled by the action of a reaction-diffusion system. In more details, we consider a classic prototype model for virus dynamics in the spatially homogeneous case derived within a framework of population dynamics~\cite{[BMSN97],[NB96]}, which is also known by the acronym SIR~\cite{[NM00],[PERT07]}. We briefly present the model which has been analytically and computationally studied in~\cite{[BPTW19]}, where a SIR type model is coupled with a Keller-Segel model.

Let us consider a May-Nowak type model which describes, by a system of ODEs, the dynamics of three components, i.e. the densities of healthy uninfected immune cells $u= u(t)$, infected immune cells $v = v(t)$, and virus particles $w = w(t)$~\cite{[NM00]}.  The model considers the following heuristic assumptions: Healthy cells are constantly produced by the body at rate $r$, die at rate $d_1 u$ and become infected on contact with the virus, at rate $\beta u w$; Infected cells are  produced at rate $\beta u w$ and die at rate $d_2 v$; New virus particles are produced at rate $kv$ and die at rate $d_3 w$. These assumptions lead to the following system of ODEs:
\begin{equation}\label{ODE}
\begin{cases}
  \displaystyle \frac{du}{dt} = - d_1\, u - \beta \, u w  + r, \hskip.5cm & t>0, \\[3mm]
     \displaystyle  \frac{dv}{dt} = - d_2\, v+ \beta \, u w,  \hskip.5cm  & t>0, \\[3mm]
      \displaystyle  \frac{dw}{dt} =- d_3\, w +k \,v,   \hskip.5cm        & t>0.
    \end{cases}
\end{equation}

This model has been quite comprehensively understood via a thorough qualitative analysis of
corresponding initial value problems (for instance cf.~\cite{[BMSN97],[NM00]}). As it is known, in addition to the  infection-free
equilibrium $Q_0 :=(\frac{r}{d_1}, 0, 0)$, if the so-called \textit{basic reproduction number} $R_0$ is greater than 1, namely  $R_0 > 1$ with
$$
R_0 := \frac{\beta kr}{d_1 \, d_2 \, d_3}\va
$$
then, the system shows an additional equilibrium  $Q^* :=(u^*,v^*, w^*)$, where
\begin{equation}
 u^*:=\frac{r}{d_1}\frac{1}{R_0}, \hskip.5cm
    v^*:=\frac{d_1 d_3}{\beta k}(R_0-1) \qquad \mbox{and}  \hskip.5cm
    w^*:=\frac{d_1}{\beta}(R_0-1).
 \end{equation}
This equilibrium  is globally asymptotically stable and positive defined, whereas if $R_0 \le 1$
then the infection-free equilibrium $Q_0$ enjoys this property~\cite{[KOR04]}.

The space dynamics in the model studied in~\cite{[BPTW19],[BT20],[FUEST19]}is modeled by a deterministic reaction diffusion dynamics acting on  $u = u(t,\bx)$, $v = v(t,\bx)$ and $w = w(t,\bx)$ which now include  space dependence:
\begin{equation}\label{1}
\begin{cases}
  \displaystyle \frac{\partial u}{\partial t} = D_u \Delta u - \chi \nabla\cdot (u\nabla v) - d_1 \, u - \beta \, u \, w + r(t, \bx), \\[3mm]
  \displaystyle \frac{\partial v}{\partial t}  =D_v \Delta v - d_2 \,v + \beta \,u \, w,   \\[3mm]
    \displaystyle \frac{\partial w}{\partial t} = D_w \Delta w - d_3\, w + k\, v,  \\
     \end{cases}
    \end{equation}
where $D_u$, $D_v$ and $D_w$ denote the respective, positive defined, diffusion coefficients and where $\chi$ represents strength and direction of the cross-diffusive interaction, while the parameters $\beta, k, d_1, d_2$, $d_3$ have been already defined above.  The reaction-diffusion action term corresponds to a simplified  Keller-Segel chemotaxis system~\cite{[BPTW19]}.

This model can be viewed as a specific example of interaction between a dynamical system modeled by ODEs and a reaction-diffusion system which creates pattern formation. Further developments may focus on the modeling of the virus dynamics, that might go beyond the limited validity of SIR models, as well as on the modeling space dynamics by selecting cross diffusion-reaction models consistent with the specific biological and physical environment
where the dynamics develops.

\section{On the micro-macro derivation of virus SIR models in a KS-system}
\label{Sec.3}

This section presents the micro-macro derivation by an asymptotic expansion somehow inspired to the Hilbert~\cite{[BC19]}.  Firstly, we derive a general kinetic model for three interacting population corresponding to the specific case of the SIR model. The derivation of macroscopic equations is treated in the following two subsections.  Firstly, we derive a general macroscopic model and  then, we show how it can be specifically referred to the virus model under consideration.

\subsection{On the derivation of a general kinetic model} \label{Subsec.(3.1)}

This subsection presents the derivation of macroscopic models, by micro-macro decomposition, of linear transport models a binary mixture of self-propelled particles  whose state, called {\sl microscopic state}, is denoted by the variable $(x, v)$,  where $x$  and $v$ are, respectively, position and velocity.  The collective description  of a mixture of particles  can be encoded in the statistical distribution functions $f_i = f_i(t, x, v)$, for $i= 1, 2, 3$. Weighted moments provide, under suitable integrability properties, the calculation of macroscopic variables.

Let us now consider the following class of equations:
\begin{equation}\label{LinearT}
\begin{cases}
 \big(\partial_t + v \cdot \nabla_{x} \big) f_1 =\nu_1\, \mathcal{T}_1[f_2,f_3](f_1)+\mu_1\,G_1(f_1,f_2,f_3,v),\\[3mm]
   \big(\partial_t + v \cdot \nabla_{x}\big)f_2=\nu_2\, \mathcal{T}_2(f_2 )+ \mu_2\,G_2(f_1,f_2,f_3,v), \\[3mm]
   \big(\partial_t + v \cdot \nabla_{x} \big) f_3= \nu_3\,  \mathcal{T}_3(f_3 )+\mu_3\,G_3(f_1,f_2,f_3,v),     \\
\end{cases}
\end{equation}
where  $G_1, G_2 , G_3 $ are interactions terms assumed depending on the quantities
$f_1, f_2, f_3$, while the operator $\mathcal{T}_i(f)$ models the dynamics of biological organisms by a velocity-jump process:
\begin{equation}\label{de}
\mathcal{T}_i(f)= \int_{V} \bigg[T_i(v^*, v)f(t, x, v^{*}) - T_i(v , v^*)f(t, x, v) \bigg]\,dv^{*}, \quad i=1,2,
\end{equation}
where $T_i(v, v^{*})$ is the probability kernel for the new velocity $v \in V$ assuming that the previous velocity was $v^{*}$.

The  derivation of macroscopic models from the kinetic model (\ref{LinearT}), can be obtained in the regime $\nu_1, \nu_2, \nu_3, \rightarrow +\infty$ corresponding to the distance between particles tending to zero. After a dimensionless of the system is obtained, see \cite{[BBNS12]}, a small parameter $\varepsilon$ can be chosen, for the parabolic scaling,  such that
\[t\longrightarrow \varepsilon t, \quad
\mu_i=\varepsilon, \quad
\nu_i=\frac{1}{\varepsilon^{q_i}},  \quad  q_i \geq 1, \quad i=1,2,3 .
\]

Then, the model (\ref{LinearT}) can be rewritten as  follows:
\begin{equation}\label{mo15}
\left\{
\begin{array}{l} \big(\varepsilon \partial_t + v \cdot \nabla_{x} \big) f_1^\varepsilon  =
\frac{1}{\varepsilon^ {q_1} } \mathcal{T}_1[f_2^\varepsilon,f_3^\varepsilon](f_1^\varepsilon )+\varepsilon\,G_1(f^\varepsilon_1,f^\varepsilon_2,f^\varepsilon_3,v),     \\[4mm]
\big(\varepsilon \partial_t + v \cdot
\nabla_{x}\big)f_2^\varepsilon =
\frac{1}{\varepsilon^{ q_2} }\mathcal{T}_2(f_2^\varepsilon )+  \varepsilon\,G_2(f^\varepsilon_1,f^\varepsilon_2,f^\varepsilon_3,v), \\[4mm]
   \big(\varepsilon \partial_t + v \cdot \nabla_{x} \big) f_3^\varepsilon  =
\frac{1}{\varepsilon^{ q_3}} \mathcal{T}_3(f_3^\varepsilon )+\varepsilon\,G_3(f^\varepsilon_1,f^\varepsilon_2,f^\varepsilon_3,v),

\end{array}\right.
\end{equation}

\vskip.2cm \noindent {\bf Assumption 3.1.} \label{assH01} The   turning operators  $\mathcal{T}_1, $ $\mathcal{T}_2, $ $\mathcal{T}_3 $  are supposed to be decomposable as
follows:
\begin{equation}
\mathcal{T}_1[f_2^\varepsilon,f_3^\varepsilon](g)= \mathcal{T}_1^{0}(g)+\varepsilon^p \, \mathcal{T}_1^{1}[f_2^\varepsilon,f_3^\varepsilon](g), \quad  p \geq 1,\label{L1}
\end{equation}
where $\mathcal{T}_1^{j}$ for $j=0,1$, is given by
\begin{equation}
\mathcal{T}_1^{j}(g)= \int_{V} \bigg[{T_1^{j}}^{*}g(t, x, v^{*}) - {T_1^{j}} g(t, x, v) \bigg]dv^{*},  \label{3.4}
\end{equation}
with ${T_1^{j}}^{*}=T_1^{j}(v^{*},v)$ and where the dependence on  $f_2$, $f_3$ , of the operator $\mathcal{T}_1 $ stems from $T_1^1 $, while we suppose that $\mathcal{T}_1^0$ is independent of $f_2$, $f_3$ ,  and  $\mathcal{T}_l$ for ($l=2,3$), is given by

\begin{equation}
\mathcal{T}_l(g)= \int_{V} \bigg[{T_l}^{*}g(t, x, v^{*}) - {T_l} g(t, x, v) \bigg]dv^{*}.
\end{equation}

\vskip.2cm \noindent {\bf Assumption 3.2.}  We
assume that the turning operators $\mathcal{T}_i(i=1,2,3)$  satisfy the following equality:
\begin{equation}\label{assH3} \int_V \mathcal{T}_i(g)dv= \int_V\mathcal{T}_1^0(g)dv=\int_V \mathcal{T}_1^1[f_2^\varepsilon,f_3^\varepsilon](g)dv=0,
\end{equation}

\vskip.2cm \noindent {\bf Assumption 3.3}\label{assH4} There exists a bounded velocity distribution $M_j(v)>0$(j=2,3) and  $M_1(v)>0$ , independent of $t, x$, such that the detailed balance
\begin{equation}
T_1^0 (v,v^{*} ) M_1(v^{*}) = T_1^0 (v^{*},v ) M_1(v),
\end{equation}
and
\begin{equation}\label{tx1}
 T_j(v,v^{*} ) M_j(v^{*}) = T_j (v^{*},v ) M_j(v),
\end{equation}
hold true. Moreover, the flow produced by these equilibrium distributions vanishes, and $M_i$ are normalized
\begin{equation}\label{TZ}
\int_V v \, M_i(v)dv  =0, \quad \int_V
M_i(v)dv =1\quad i=1,2,3.
\end{equation}
In addition, we assume that the kernels  $T_j(v,v^{*})$ and  $T_1^0(v,v^{*})$  are bounded and that there exist constants  $\sigma_j>0$ and $\sigma_1>0$, $j=2,3$, such that
\begin{eqnarray}\label{cx}
T_j(v,v^{*})\geq \sigma_jM_j(v),\quad T_1^0(v,v^{*})\geq \sigma_1 M_1(v),
\end{eqnarray}
for all $ (v,v^{*}) \in V\times V $, $ x\in \Omega$ and $ t>0$.
\vskip.2cm

Given that $L_j=\mathcal{T}_j(j=2,3)$ and $L_1=\mathcal{T}^0_1$. Technical calculations yields the following Lemma:
\vskip.1cm
\begin{Lemma}\label{LE1}
Suppose that  Assumptions 3.3 holds. Then, for $i=1,2,3,$ the following properties of the operators $L_1$, $L_2$ and $L_3$  hold:
\begin{itemize}
\item[i)] The operator $L_i$ is self-adjoint in the space $\displaystyle{{L^{2}\left(V ,{dv \over M_i}\right)}}$.
\vskip.1cm \item[ii)]
For $f\in L^2$, the equation $L_i(g)=f$ has a unique solution $\displaystyle{g \in L^{2}\left(V, {dv \over M_i}\right)}$, which satisfies
$$
 \int_{V} g(v)\, dv = 0 \quad  \hbox{if and only if} \quad   \int_{V} f(v)\, dv =0.
$$
\vskip.1cm \item[iii)]  The equation $L_i(g) =v \,  M_i(v)$ has a unique solution that we call $\theta_i(v)$.
\vskip.1cm \item[iv)]  The kernel of $L_i$ is $N(L_i) = vect(M_i(v))$.
\end{itemize}
\end{Lemma}

\subsection{Derivation of a general  macroscopic models} \label{Subsec.(3.2)}

A system coupling a hydrodynamic part with a kinetic part of the distribution functions, is derived in this subsection. Then it is proved that  such a system is equivalent to the two scale kinetic equation (\ref{mo15}). This new formulation provides the basis for the derivation of the  general model we are looking for.

In the remainder, the integral with respect to the variable $v$ will be denoted by $\langle \cdot  \rangle$. This notation is used also for  any argument  within $\langle \, \rangle$. In addition,   let us denote by  $f=(f_1, f_2,f_3)$ the solution of (\ref{mo15}), where $f$  is decomposed as follows:
\begin{equation} \label{eq}
f_{1}^\varepsilon(t,x,v)=\sum^{q_1+1}_{i=0}\varepsilon^{i}g_{i}(t,x,v)+O(\varepsilon^{q_1+2}),
  \end{equation}
\begin{equation} \label{eq1}
f_{2}^\varepsilon(t,x,v)=\sum^{q_1+1}_{j=0}\varepsilon^{j}h_{j}(t,x,v)+O(\varepsilon^{q_2+2}).
 \end{equation}
and
\begin{equation} \label{eq2}
f_{3}^\varepsilon(t,x,v)=\sum^{q_3+1}_{l=0}\varepsilon^{l}k_{l}(t,x,v)+O(\varepsilon^{q_3+2}).
 \end{equation}

In order to develop asymptotic analysis of Eq.~(\ref{mo15}),  additional assumptions on the operator $\mathcal{T}_1^1$  and the interaction terms
$G_{\hat{i}} (\hat{i}=1,2,3)$   are needed.

\vskip.2cm \noindent {\bf Assumption 3.4.}
We assume that the turning operator $\mathcal{T}_1^1$ and the interaction terms
$G_{\hat{i}}(i=1,2,3)$ satisfy the following asymptotic behavior as:
\begin{equation}
\begin{cases}
\mathcal{T}_1^1\big[f_{2}^\varepsilon, f_{3}^\varepsilon](g)= \mathcal{T}_1^1[h_0,k_0] (g)\\[2mm]
\hskip1cm  + \displaystyle{\sum^{q_1+1-p}_{m=1}\varepsilon^m \mathcal{R}_1^m [h_0,....,h_j,k_0,...,,k_l](g)}\\[2mm]
\hskip1cm + O(\varepsilon^{q_1+2-p}) ,~~ \forall p\leq q_1,~~ \forall g,h_j,k_l,\\[4mm]
\mathcal{T}_1^1\big[f_{2}^\varepsilon, f_{3}^\varepsilon](g)= \mathcal{T}_1^1[h_0,k_0](g)+ O(\varepsilon) ,\quad \forall p> q_1,~~ \forall g,\\
\end{cases}
\end{equation}
for $j=0,1,...,q_2+1$ and $l=0,1,...,q_3+1$.

\vskip.2cm

Then,  (\ref{assH3}) rapidly yields:
\begin{equation}\label{R1}
\int_{V}  \mathcal{R}_1^m[h_0,.....,h_j,k_0,....,,k_l](\varphi)dv=0,  \quad \forall \varphi.
 \end{equation}
and
\begin{equation}\label{G2}
 G_{\hat{i}}(g+\varepsilon \,  \hat{g}, h+\varepsilon  \, \hat{h}, k+\varepsilon  \, \hat{k},v) =G_{\hat{i}}(g, h,k,v)+ O(\varepsilon),\quad\forall g,\hat{g},h,\hat{h},k,\hat{k}
 \end{equation}
for ${\hat{i}}=1,2,3$.

Then, the first terms of Hilbert expansion of equal order in $\varepsilon^i$, $\varepsilon^j$   and  $\varepsilon^l$  for $i=0,1,...,q_1+1$, $j=0,1,...,q_2+1$ and $l=0,1,...,q_3+1$ are:
\begin{equation}\label{c1}
\ve^0:
 \left\{
 \begin{array}{l}
\mathcal T_{1}^{0}(g_{0})=0,\\[4mm]
\mathcal T_{2}(h_{0})=0,\\[4mm]
\mathcal T_{3}(k_{0})=0,\\
\end{array}  \right.
\end{equation}
\begin{equation}\label{ccoo}
\ve^1:
 \left\{
 \begin{array}{l}
\mathcal T_{1}^{0}(g_{1})=\delta_{q_1,1}\, v\cdot \nabla_{x} g_{0}-\delta_{p,1}\mathcal T_{1}^{1}[h_{0},k_{0}](g_{0}),\\[4mm]
\mathcal T_{2}(h_{1})=\delta_{q_2,1} v\cdot \nabla_{x} h_{0},\\[4mm]
\mathcal T_{3}(k_{1})=\delta_{q_3,1} v\cdot \nabla_{x} k_{0},\\
\end{array}  \right.
\end{equation}
\begin{equation}\label{ccc}
{\ve^2:}
 \left\{
 \begin{array}{llll}
 \mathcal T_{1}^{0}(g_{2})=\delta_{q_1,1} \big( \partial_t g_{0}+ v\cdot \nabla_{x} g_{1}\big)-\delta_{p,2} \mathcal T_{1}^{1}[h_{0}, k_{0}](g_{0})-\delta_{p,1} \mathcal T_{1}^{1}[h_{0}, k_{0}](g_{1})\\[3mm]
\hskip1cm +\delta_{q_1,2} v\cdot \nabla_{x} g_{0}    -\delta_{p,1} \mathcal R_{1}^{1}[h_{0},...,h_{j},k_{0},...,k_{l}](g_{0})\\[3mm]
 \hskip2cm -\delta_{q_{1},1}G_{1}(g_{0},h_{0},k_{0} ,v), \quad \forall p\leq q_1,\\[4mm]
 \mathcal T_{1}^{0}(g_{2})=\delta_{q_1,1} \big( \partial_t g_{0}+ v\cdot \nabla_{x} g_{1}\big)-\delta_{p,2} \mathcal T_{1}^{1}[h_{0}, k_{0}](g_{0})\\[3mm]
\hskip1cm +\delta_{p,2} v\cdot \nabla_{x} g_{0} - \delta_{q_{1},1}G_{1}(g_{0},h_{0},k_{0} ,v), \quad \forall p>q_1,\\[4mm]
 \mathcal T_{2}(h_{2})= \delta_{q_2,1} \big(\partial_t h_{0} + v\cdot \nabla_{x} h_1 \big) + \delta_{q_2,2}\, v\cdot \nabla_{x}{h_0}\\[3mm]
 \hskip2cm -\delta_{q_2,1}G_{2}(g_{0},h_{0},k_{0}), \\[4mm]
 \mathcal T_{3}(k_{2})=
\delta_{p,1} \big(\partial_t k_{0} + v\cdot \nabla_{x} k_1 \big) + \delta_{q_3,2}\, v\cdot \nabla_{x}{k_0} \\[3mm]
 \hskip2cm -\delta_{q_3,1}G_{3}(g_{0},h_{0},k_0).
\end{array}
 \right.
\end{equation}

Further calculations yield:
\begin{equation}\label{ccc}
{\varepsilon^{q_{1}+1}:}
 \left\{
 \begin{array}{llll}
 \mathcal T_{1}^{0}(g_{q_{1}+1})=\partial_t g_{0}+ v\cdot \nabla_{x} g_{1}-
\sum^{q_1+1-p}_{m=1}\bigg[\delta_{p,q_1+1-m} \mathcal T_{1}^{1}[h_{0}, k_{0}](g_{m})\\[2mm]
\hskip1cm -  \sum^{q_1-p-m}_{i=1}  \delta_{p,q_1-m-i} \mathcal R_{1}^{i}[h_{0},...,h_{j},k_{0},...,k_{l}](g_{m})\bigg]\\[2mm]
\hskip1cm - G_{1}(g_{0},h_{0},k_{0} ,v)\\[2mm]
 \hskip1cm  -  \sum^{q_1+1-p}_{i=1}  \delta_{p,q_1+1-i}\mathcal R_{1}^{i}[h_{0},...,h_{j},k_{0},...,k_{l}](g_{0}),\\[3mm]
 \hskip2cm \forall p\leq q_1,\\[4mm]
 \mathcal T_{1}^{0}(g_{q_{1}+1})=\partial_t g_{0}+ v\cdot \nabla_{x} g_{1} - G_{1}(g_{0},h_{0},k_{0} ,v)\\[2mm]
\hskip1cm -\delta_{p,q_1+1} \mathcal T_{1}^{1}[h_{0}, k_{0}](g_{0}),\\
\end{array}
 \right.
\end{equation}
and
\begin{equation}\label{ccc}
 \left\{
 \begin{array}{llll}
\varepsilon^{q_{2}+1}:\mathcal T_{2}(h_{q_{2}+1})&=&\partial_t h_{0}+ v\cdot \nabla_{x} h_{1}-G_{2}(g_{0},h_{0},k_{0} ,v),\\[4mm]
\varepsilon^{q_{3}+1}:\mathcal T_{3}(k_{q_{3}+1})&=&\partial_t k_{0}+ v\cdot \nabla_{x} k_{1}-G_{3}(g_{0},h_{0},k_{0} ,v),\\
\end{array}
 \right.
\end{equation}
where $\delta_{a,b}$ stands for the Kronecker delta.

The first equation of (\ref{c1}) implies that
$$
g_{0}\in  vect(M_{1}(v)),  h_{0}\in  vect(M_{2}(v)), \hskip1cm  \hbox{and} \hskip1cm k_{0}\in  vect(M_{3}(v)).
$$
Therefore
$\exists c(t,x), \exists s(t,x), \exists u(t,x)$ such that
\begin{equation} \label{nS1}   g_0(t,x,v)= M_1(v)\,c(t,x),
\end{equation}
\begin{equation}
\label{nS12}h_0(t,x,v)= M_2(v)\,s(t,x)~~\hbox{and} ~~k_0(t,x,v)= M_3(v)\,u(t,x).
\end{equation}

Using (\ref{assH3}), (\ref{TZ}) and  (\ref{nS1})-(\ref{nS12}), we conclude that Eq. (\ref{ccoo}) satisfies the solvability condition, therefore $g_{1}$ and $h_{1}$ are given by
\begin{equation}\label{ph}
  \begin{cases}
 g_{1}=\delta_{q_1,1} (\mathcal T_{1}^{0})^{-1}(v\cdot\nabla_{x} g_{0})- \delta_{p,1}(\mathcal T_{1}^{0})^{-1}(\mathcal T_{1}^{1}[h_{0}](g_{0})),
 \\[3mm]
 h_{1}=\delta_{q_2,1} \mathcal T_{2}^{-1}(v\cdot\nabla_{x} h_{0}),
 \\[3mm]
 k_{1}=\delta_{g_3,1} \mathcal T_{3}^{-1}(v\cdot\nabla_{x} k_{0}).
\end{cases}
\end{equation}

The calculations of $ g_{q_1+1}$, $ h_{q_2+1}$, and  $ k_{q_3+1}$  are obtained from the solvability conditions at $O(\varepsilon^{q_1+1})$, $O(\varepsilon^{q_2+1})$  and $O(\varepsilon^{q_3+1})$,  which are given by the following:
\begin{equation}\label{cc1}
 \left\{
 \begin{array}{llll}
\displaystyle\int_{V}\bigg(\partial_t g_{0}+ v\cdot \nabla_{x} g_{1}-
\sum^{q_1+1-p}_{m=1}\bigg[\delta_{p,q_1+1-m} \mathcal T_{1}^{1}[h_{0}, k_{0}](g_{m})\\[3mm]
-  \sum^{q_1-p-m}_{i=1}  \delta_{p,q_1-m-i} \mathcal R_{1}^{i}[h_{0},...,h_{j},k_{0},...,k_{l}](g_{m})\bigg]-G_{1}(g_{0},h_{0},k_{0} ,v)\\[3mm]
-  \sum^{q_1+1-p}_{i=1}  \delta_{p,q_1+1-i}\mathcal R_{1}^{i}[h_{0},...,h_{j},k_{0},...,k_{l}](g_{0})\bigg)dv=0,\\[3mm]
\hskip3cm \forall p\leq q_1,\\[4mm]
{}
\displaystyle\int_{V}\bigg( \partial_t g_{0}+ v\cdot \nabla_{x} g_{1}-G_{1}(g_{0},h_{0},k_{0} ,v)-\delta_{p,q_1+1} \mathcal T_{1}^{1}[h_{0}, k_{0}](g_{0})\bigg)dv=0,\\[3mm]
\hskip3cm \forall p>q_1,
\end{array}
 \right.
\end{equation}
and
\begin{equation}\label{cc2}
 \left\{
 \begin{array}{llll}
\displaystyle\int_{V}\bigg(\partial_t h_{0}+ v\cdot \nabla_{x} h_{1}-G_{2}(g_{0},h_{0},k_{0} ,v)\bigg)dv=0,\\[4mm]
\displaystyle\int_{V}\bigg(\partial_t k_{0}+ v\cdot \nabla_{x} k_{1}-G_{3}(g_{0},h_{0},k_{0} ,v)\bigg)dv=0,
\end{array}
 \right.
\end{equation}

Using (\ref{TZ}), (\ref{assH3}), (\ref{R1})  and (\ref{nS1})-(\ref{ph}), denoting by $<\cdot>$ the integral with respect to the variables $v$, shows that the system (\ref{cc1})-(\ref{cc2})  can be  rewritten  as follows:
\begin{equation}\label{aa}
\left\{
 \begin{array}{llll}
 \displaystyle
\partial_t c +\delta_{q_1,1}\, \left\langle v \cdot \nabla_{x}(\mathcal T_{1}^{0})^{-1}(v M_{1}\cdot\nabla_{x} c)\right\rangle \\[4mm]
 \hskip2cm -\delta_{p,1}\,\left\langle (\mathcal T_{1}^{0})^{-1}(\mathcal T_{1}^{1}[M_{2} s,M_{3} u](M_{1} c))\right\rangle\\[4mm]
\hskip2cm - \left\langle G_{1}(M_{1} c,M_{2} s,M_{3} u,v)\right\rangle = 0,  \vspace{.2cm}\\[4mm]
	\displaystyle \partial_t s+
\delta_{q_2,1}\left\langle v\cdot\nabla_x\mathcal T_{2}^{-1}(v M_{2}\cdot\nabla_{x} s)\right\rangle  -  \left\langle G_{2}(M_{1} c,M_{2} s,M_{3} u,v)\right\rangle \vspace{.2cm} = 0,\\[4mm]
	\displaystyle \partial_t u+
\delta_{q_3,1}\left\langle v\cdot\nabla_x\mathcal T_{3}^{-1}(v M_{3}\cdot\nabla_{x} u)\right\rangle  -  \left\langle G_{3}(M_{1} c,M_{2} s,M_{3} u,v)\right\rangle \vspace{.2cm} = 0.
\end{array} \right.
\end{equation}

As $\mathcal{T}_{1}^{0},$ $\mathcal{T}_2$ and $\mathcal{T}_3$  are self-adjoint operators in $L^{2}\left(D_v ,{dv \over M_1(v)}\right),$ $L^{2}\left(D_v ,{dv \over M_2(v)}\right)$ and
$L^{2}\left( D_v ,{dv\, \over M_3(v)}\right)$,  one has the following computations:
$$
\left\langle v. \nabla_{x} (\mathcal{T}_{1}^{0})^{-1}(v M_1 \cdot \nabla_{x} c)\right\rangle=\div_x\bigg( \langle v \otimes
\theta_1(v)\rangle\cdot \nabla_{x} c \bigg),
$$
$$
\left\langle v\cdot \nabla_{x} \mathcal{T}_2^{-1}(v M_2 \cdot \nabla_{x} s)\right\rangle= \div_x\bigg( \langle v \otimes
\theta_2(v)\rangle\cdot \nabla_{x} s \bigg),
$$
$$
\left\langle v\cdot \nabla_{x} \mathcal{T}_3^{-1}(v M_3 \cdot \nabla_{x} u)\right\rangle= \div_x\bigg( \langle v \otimes
\theta_3(v)\rangle\cdot \nabla_{x}u \bigg),
$$
and
$$
 \left\langle v\cdot \nabla_{x}
(\mathcal{T}_{1}^{0})^{-1}(\mathcal{T}_1^1[M_2s,M_3u](M_1c))\right\rangle
 = \div_x \left\langle \frac{\theta_1(v)}{M_1(v)}c\mathcal{T}_1^1[M_2s,M_3u](M_1)\right \rangle,
 $$
 where $\theta_1$ and $\theta_2$ are given in Lemma 2.

Therefore, the macroscopic model (\ref{aa}) can be written as follows:
\begin{equation}\label{mM22}
\left\{
\begin{array}{l}
 \partial_t c +  \div_{x} \, (\delta_{p,1}\,  c\, \alpha(s,u)- \delta_{q_1,1}\, D_c \cdot \nabla_x c)-  H_1(c,s,u) =0,  \\[4mm]
    \partial_t s-  \delta_{q_2,1}\, \div_{x} \, (  D_s \cdot \nabla_x s)-  H_2(c,s,u) =0,  \\[4mm]
    \partial_t u  - \delta_{q_3,1} \,\div_{x} \, (    D_u\cdot \nabla_x u)-  H_3(c,s,u) =0,  \\[4mm]
\end{array} \right. \end{equation}
where  $ D_c$, $D_s$, $ D_u$, $\alpha$ are given, respectively, by
\begin{equation}\label{di}
D_c =- \int_V v \otimes \theta_1(v) dv, \quad D_s =- \int_V v \otimes \theta_2(v) dv,  \quad D_u =- \int_V v \otimes \theta_3(v) dv,
\end{equation}
and
\begin{equation} \label{w1}
\alpha(s,u)= - \int_V  {\theta_1(v)\over M_1(v)}  T_{1}^{1}[M_{2} s,M_{3} u](M_{1} c))dv,
\end{equation}
while  $H_i(c,s,u)$, (i=1,2,3) are given by:
\begin{equation} \label{w32}
 H_i(c,s,u)= \int_V G_{i}(M_1 c, M_2 s,M_3 u,v) dv.
\end{equation}

\subsection{Derivation of virus models with in a Keller-Segel system}\label{Subsec.(3.3)}

 More in detail, let us  consider the following kernels:
\begin{equation}\label{TO}
T_1^0(v,v^{*})=\sigma_1 M_1(v), \quad T_2(v,v^{*})=\sigma_2 M_2(v), \quad T_3(v,v^{*})=\sigma_3M_3(v),
\end{equation}
with $\sigma_1, \sigma_2, \sigma_3>0$. 

Hence, the leading turning operators $\mathcal{T}_j$ (j=2,3) and  $\mathcal{T}^0_1$  can be viewed as relaxation operators:
\begin{equation}\label{relaxation1}
\mathcal{T}_1^{0}(g)= -\sigma_1 \Big(g- M_1 \langle g\rangle\Big),
 \end{equation}
 \begin{equation}\label{relaxation1b}
 \mathcal{T}_2(g)= -\sigma_2\Big(g-  M_2 \langle g\rangle\Big),
\end{equation}
\begin{equation}\label{relaxation2}
\mathcal{T}_2(g)= -\sigma_2\Big(g-  M_2 \langle g\rangle\Big),
\end{equation}
and
\begin{equation}\label{relaxation3}
\mathcal{T}_3(g)= -\sigma_3 \Big(g- M_3 \langle g\rangle\Big).
\end{equation}

Moreover, $\theta_1$, $\theta_2$, and $\theta_3$  are given by
\[
\theta_1(v)= -\frac{1}{ \sigma_1}\,  v M_1 (v)\va  \quad \theta_2(v)= -\frac{1}{ \sigma_2}\,  v M_2(v),  \quad  \hbox{and} \quad   \theta_3(v)= -\frac{1}{ \sigma_3}\,  v M_3(v)\va
\]
while $\alpha$, are defined by (\ref{w1})--(\ref{w32}), and are computed as follows:
\begin{equation}
\alpha( s,u)=  \frac{1}{\sigma_1} \int_V v \mathcal{T}_1^1[M_2 s,M_3 u](M_1(v) )dv, \label{alpha}
\end{equation}

The diffusion tensors $D_c$, $D_s$, and  $D_u$  are given by
\begin{equation}\label{df}
D_c= \frac{1}{\sigma_1}\int_V v \otimes v M_1(v) dv,  \quad \hbox{and} \quad D_s= \frac{1}{\sigma_2}\int_V v \otimes v M_2(v) dv,
 \end{equation}
\begin{equation}\label{df1}
D_u= \frac{1}{\sigma_3}\int_V v \otimes v M_3(v) dv,
 \end{equation}
while $H_1, H_2$, $H_3$  are still given by (\ref{w32}).

\vskip.2cm  \noindent $\bullet$  Let us now consider that $q_1=q_2=q_3=2$  and $p=2$, then from (\ref{mM22}) one has the following macro-scale type models up to $\varepsilon$:
\begin{equation}\label{mM501}
    \left\{ \begin{array}{ll}
  \displaystyle \frac{dc}{dt} = H_1( c, s,u),  \\[4mm]
     \displaystyle  \frac{ds}{dt} = H_2( c, s,u),  \\ [4mm]
      \displaystyle  \frac{du}{dt} =H_3( c, s,u).
    \end{array} \right.
\end{equation}

 The role of the terms $H_1(  c, s,u)$, $H_2(  c, s,u)$, and $H_3(  c, s,u)$  in
(\ref{w32}) consists in modeling the interaction between the for quantities of the mixture.  For example,  by choosing:
\begin{equation} \label{KK}
G_1(f_1, f_2,f_3,v)=-\frac{d_1}{|V|} \,\frac{f_1}{M_1} -\frac{\beta}{|V|}\frac{f_1}{M_1} \frac{f_3}{M_3}+\frac{1}{|V|}\,r,
\end{equation}

\begin{equation} \label{KK1}
G_2(f_1, f_2,f_3,v)=-\frac{d_2}{|V|} \,\frac{f_2}{M_2} +\frac{\beta}{|V|}\frac{f_1}{M_1} \frac{f_3}{M_3},
\end{equation}

\begin{equation} \label{KK2}
G_3(f_1, f_2,f_3,v)=-\frac{d_3}{|V|} \,\frac{f_3}{M_3} +\frac{k}{|V|}\frac{f_2}{M_2},
\end{equation}

Therefore, the macroscopic model (\ref{mM501}) writes:
\begin{equation}\label{ODE1}
    \left\{ \begin{array}{ll}
  \displaystyle \frac{dc}{dt} = - d_1\, c - \beta \, c  u  + r, \\[4mm]
     \displaystyle  \frac{ds}{dt} = - d_2\, s+ \beta \, c u, \\ [4mm]
      \displaystyle  \frac{du}{dt} =- d_3\,u +k \,s.
    \end{array} \right.
\end{equation}

\vskip.2cm  \noindent $\bullet$   Let us also consider that $q_1=q_2=q_3=1$, $p=1$, and  the following choice:
\begin{equation}\label{Q}
T_1^1[f_2,f_3]= K_{\frac{f_2}{M_2}}(v, v^{*})\cdot \nabla_{x} \frac{f_2}{M_2},
\end{equation}
where $K_{\frac{f_2}{M_2}}(v, v^{*})$ is a vector valued function satisfying the following:
 \begin{equation}\label{K}
K_{S+\varepsilon  \, {g \over M_2}}= K_S +O(\varepsilon), \quad \hbox{as}
\quad \varepsilon \rightarrow 0.
\end{equation}
 Then $\mathcal T_1^1$ satisfies (\ref{TZ}), and leads to the following:
\begin{equation*}
\mathcal{T}^1_1[M_2s,M_3u](M_1)= \psi(v,s)\cdot \nabla_x s,
\end{equation*}
where
\begin{equation}\label{Qj}
\psi(v,s)=\int_V \Big( K_s(v, v^{*})M_2(v^{*})  -  K_s( v^{*}, v)M_2(v)\Big)dv^{*}.
\end{equation}

Finally,  $\alpha(s)$, defined in \eqref{alpha}, is given by $\alpha(s)= \chi(s)\cdot \nabla_x s$, where the chemotactic sensitivity $\chi(s)$ is given by the matrix
\begin{equation}\label{C}
\chi(s)=\frac{1}{\sigma_1} \int_V v\otimes \psi(v, s) dv.
\end{equation}

Therefore, the macroscopic model (\ref{mM22}) can be written as follows:
\begin{equation}\label{mM5}
\begin{cases}
 \partial_t c +  \div_{x} \, \bigg(c\,  \chi(s)\cdot \nabla_x s- \,  (D_c \cdot \nabla_x c)\bigg) = H_1( c, s,u), \\[4mm]
       \partial_t s- \div_{x} \, (D_s \cdot \nabla_x s) = H_2( c, s,u), \\[4mm]
    \partial_t u -  \div_{x} \, (D_u \cdot \nabla_x u)= H_3( c, s,u).
   \end{cases}
   \end{equation}

Supposing that the $G_ i (i=1,2,3),$  are given by (\ref{KK})-(\ref{KK2}), yields the following macroscopic model (\ref{mM5}):
\begin{equation}\label{ODE1}
\begin{cases}
  \displaystyle \partial_t c + \div_{x} \, \bigg(c\,  \chi(s)\cdot \nabla_x s- \,  (D_c \cdot \nabla_x c)\bigg) = - d_1\, c - \beta \, c \, u  + r, \\[4mm]
     \displaystyle   \partial_t s - \div_{x} \,(D_s \cdot \nabla_x s)= - d_2\, s+ \beta \, c\, u, \\[4mm]
      \displaystyle   \partial_t - \div_{x} \,(D_u \cdot \nabla_x u) =- d_3\,u +k \,s,
    \end{cases}
\end{equation}
where the $D_c ,D_s , D_u$, and $\chi(s)$  are given by (\ref{df})-(\ref{df1})  and (\ref{C}).

\section{Critical analysis and perspectives }\label{Sec.}

 We trust that the micro-macro derivation developed in our paper, can be further extended to a variety of exotic models including models of the dynamics of different types of virus. This subsection simply introduces this topic which will be developed in a well defined research program. Indeed, the complexity of the dynamics may substantially increase thus requiring nontrivial developments of the approach. In more details on virus  dynamics, possible developments might account for  delay-distributed dynamics~\cite{[EA17]}, infection model with multi-target cells~\cite{[EAA19]}, stochastic models which include selective progression~\cite{[GEA17]}. In addition, the recent Covid-19 virus pandemic has generated a huge number of models that can be viewed as a technical development of SIR type models. The research article~\cite{[GBM20]} is one of the first contributions to this topic. Subsequently, it has been followed by various papers which have gone beyond the framework of compartmental models
accounting for multiscale features and of the immune competition inside the lung~\cite{[BBC20]}, transport dynamics~\cite{[BDP21]}, contagion in crowds~\cite{[KQ20]}, and various others.



\begin{thebibliography}{99}

\bibitem{[BBNS12]}
\newblock   N.~Bellomo, A.~Bellouquid, J.~Nieto, and J.~Soler,
\newblock  On the asymptotic theory from microscopic to  macroscopic tissue models: an overview with perspectives,
\newblock   \textit{Mathematical Models Methods Applied Sciences},  \textbf{22} (2012), paper n.~1130001.

\bibitem{[BBC20]}
\newblock N.~Bellomo, R.~Bingham, M.A.J.~Chaplain, G.~Dosi, G.~Forni, D.A.~Knopoff, J.~Lowengrub, R.~Twarock, and M.E.~Virgillito,
\newblock A multi-scale model of virus pandemic: Heterogeneous interactive entities in a globally connected world,
 \newblock \textit{Mathematical Models and Methods Applied Sciences}, \textbf{30} (2020), 1591--1651.


\bibitem{[BBDG21]} 
\newblock  N.~Bellomo,  D.~Burini, G.~Dosi, L.~Gibelli, D.Knopoff, P.~Terna, and M.E.~Virgillito,
\newblock What is life? A perspective of the mathematical kinetic theory of active particles,
\newblock   \textit{Mathematical Models Methods Applied Sciences},  \textbf{31} (2021),  1821--1866.


\bibitem{[BOSTW21]}
\newblock  N.~Bellomo,  N.~Outada, J.~Soler, Y.~Tao, and M.~Winkler,
\newblock  Chemotaxis and cross diffusion models in complex environments: Modeling  towards a multiscale vision,
\newblock  \textit{Mathematical Models Methods Applied Sciences}, to appear, (2022).

\bibitem{[BPTW19]}
\newblock N.~Bellomo, K.~Painter, Y.~Tao and M.~Winkler,
\newblock Occurrence vs. absence of taxis-driven instabilities in a May--Nowak model for virus infection,
\newblock  \textit{SIAM Journal Applied Mathematics}, \textbf{79(5)}  (2019),  1990--2010.



\bibitem{[BT20]}
  \newblock N.~Bellomo and Y.~Tao,
  \newblock Stabilization in a chemotaxis model for virus infection,
  \newblock \textit{Discrete Continuous Dynamical Systems Series~S}, \textbf{13}  (2020), 105--117.

\bibitem{[BMSN97]}
  \newblock  S.~Bonhoeffer, R.M.~May, G.M.~Shaw, and M.A.~Nowak,
  \newblock Virus dynamics and drug therapy,
  \newblock \textit{Proceedings National Academy Sciences~USA}, \textbf{94}  (1997), 6971--6976.

\bibitem{[BDP21]}
\newblock W.~Boscheri, G.~Dimarco, and L.~Pareschi,
\newblock Modeling and simulating the spatial spread of an epidemic through multiscale kinetic transport equations,
\newblock \textit{Mathematical  Models and Methods in Applied Sciences},  \textbf{31},  (2021).
https://doi.org/10.1142/S0218202521400017.


\bibitem{[BC17]}
\newblock D.~Burini and N.~Chouhad,
\newblock Hilbert method toward a multiscale analysis from kinetic to macroscopic models for active particles,
\newblock \textit{Mathematical Models Methods Applied Sciences}, \textbf{27(7)}  (2017), 1327--1353.

\bibitem{[BC19]}
\newblock D.~Burini and N.~Chouhad,
\newblock A Multiscale view of nonlinear diffusion in biology: From cells to tissues,
\newblock \textit{Mathematical Models Methods Applied Sciences}, \textbf{29(4)} (2019), 791--823.

\bibitem{[EA17]}
\newblock A.M.~Elaiw and N.H.~AlShamrani,
\newblock Stability of a general delay-distributed virus dynamics model with multi-staged infected progression and immune response,
\newblock \textit{Mathematical Methods Applied Sciences}, \textbf{40} (2017), 699--719.


\bibitem{[EAA19]}
\newblock A.M.~Elaiw,  T.O.~Alade and S.M.~Alsulami,
\newblock Global dynamics of delayed CHIKV infection model with multitarget cells.
\newblock \textit{Journal Applied Mathematics and Computing}, \textbf{60} (2019), 303--325.


\bibitem{[FUEST19]}
\newblock M.~Fuest,
\newblock Boundedness enforced by mildly saturated conversion in a chemotaxis--May--Nowak model for virus infection,
\newblock \textit{Journal of Mathematical Analysis and Applications}, \textbf{472(2)} (2019), 1729--1740.

\bibitem{[GBM20]}
\newblock M.~Gatto, E.~Bertuzzo, L.~Mari, S.~Miccoli, L.~Carraro, R.~Casagrandi, and A.~Rinaldo,
\newblock Spread and dynamics of the COVID-19 epidemic in Italy: Effects of emergency containment measures,
\newblock \textit{Proceedings of the National Academy of Sciences}, \textbf{117(19)}, 10484--10491, (2020).


\bibitem{[GEA17]}
\newblock L.~Gibelli,  A.M.~Elaiw and   M.A.~Alghamdi,
\newblock Heterogeneous population dynamics of active particles: Progression, mutations, and selection dynamics,
\newblock \textit{Mathematical Models Methods Applied Sciences}, \textbf{27} (2017), 617--640.


\bibitem{[GK14]}
\newblock  A.N.~Gorban and I.~Karlin,
\newblock  Hilbert's $6$th problem: exact and approximate hydrodynamic manifolds for kinetic equations,
\newblock \textit{Bulletin American Mathematical Society}, \textbf{51} (2014), 187--246.

\bibitem{[HILBERT]}
\newblock D.~Hilbert,
\newblock Mathematical problems,
\newblock \textit{Bulletin American Mathematical Society}, \textbf{8(10)} (1902), 437--479.

\bibitem{[HP09]}
\newblock T. Hillen and K.J. Painter,
\newblock A user's guide to PDE models for chemotaxis,
\newblock \textit{Journal Mathematical Biology},  \textbf{58} (2009), 183--217.


\bibitem{[KS70]}
\newblock  E.F.~Keller and L.A.~Segel,
\newblock Initiation of slime mold aggregation viewed as an instability,
\newblock \textit{Journal  Theoretical Biology}, \textbf{26} (1970), 399--415.

\bibitem{[KS71]}
\newblock E.F.~Keller and L.A.~Segel,
\newblock  Model for chemotaxis,
 \newblock \textit{Journal  Theoretical Biology}, \textbf{30} (1971), 225--234.

\bibitem{[KQ20]}
D.~Kim and A.~Quaini,
Coupling kinetic theory approaches for pedestrian dynamics and disease contagion in a confined environment.
\textit{Math. Mod. Meth. Appl. Sci.} \textbf{30(10)} (2020), 1893--1915.


 \bibitem{[KOR04]}
  \newblock A.~Korobeinikov,
  \newblock Global properties of basic virus dynamics models,
  \newblock  \textit{Bulletin Mathematical Biology}, \textbf{66} (2004),  879--883.



\bibitem{[NB96]}
  \newblock M.A.~Nowak and C.R.M.~Bangham,
  \newblock Population dynamics of immune responses to persistent viruses,
  \newblock   \textit{Science}, \textbf{272} (1996), 74--79.

\bibitem{[NM00]}
  \newblock M.A.~Nowak and  R.~May,
  \newblock \textbf{Virus dynamics: Mathematical principles of immunology and virology},
  \newblock  Oxford University Press, (2000).



\bibitem{[PERT07]}
\newblock B.~Perthame,
\newblock \textbf{Transport Equations in Biology},
\newblock Birkh\"auser Basel, (2007).




 \end{thebibliography}
\end{document}